# Structural-equation Estimation without Instrumental Variables

Eric Blankmeyer

Email eb01@txstate.edu

September 2017

**Abstract**. For a single equation in a system of linear equations, estimation by instrumental variables is the standard approach. In practice, however, it is often difficult to find valid instruments. This paper proposes a maximum-likelihood method that does not require instrumental variables; it is illustrated by simulation and with a real data set.

**Introduction**

The consistent estimation of a simultaneous equation is a significant topic in econometrics and a challenging problem at the practical level. For an equation embedded in a system of simultaneous linear equations, estimation by ordinary least squares (OLS) is generally inconsistent; and the apparently-universal recommendation is to use instrumental variables (IV), for example two-stage least squares (TSLS) or limited information maximum likelihood (LIML). "The substantial literature on estimation with weak instruments has not yet produced a serious practical competitor to the usual instrumental variable estimator" (Leamer 2010, 44-45). "The various methods that have been developed for simultaneous-equations models are all IV estimators" (Greene 2003, 398).

However, the IV approach is frequently problematic since a valid instrument must be uncorrelated with the random component of the variable for which it is a proxy but sufficiently correlated with its systematic component.

> Those who use instrumental variables would do well to anticipate the inevitable barrage of questions about the appropriateness of their instruments. Ever-present asymptotic bias casts a large dark shadow on instrumental variables estimates and is what limits the applicability of the estimate even to the setting that is observed, not to mention extrapolation to new settings. In addition, small sample bias of instrumental variables estimators, even in the consistent case, is a huge neglected problem with practice…(Leamer 2010, 35)

> What can go wrong with instrumental variables? The most important potential problem is a bad instrument, that is, an instrument that is correlated with the omitted variables (or the error term in the structural equation of interest in the case of simultaneous equations). Especially worrisome is the possibility that an association between the instrumental variable and omitted variables can lead to a bias in the resulting estimates that is much greater than the bias in ordinary least squares estimates.

> Moreover, seemingly appropriate instruments can turn out to be correlated with omitted variables on closer examination (Angrist and Krueger 2001, 79).

Widespread concern about the effectiveness of IV methods is especially apparent from the voluminous literature on weak instruments that has emerged during the last twenty years, a small but distinguished sample of which includes Staiger and Stock (1997), Angrist and Krueger (2001), Mariano (2001), Stock et al. (2002), Andrews and Stock (2006), Chernozhukov and Hansen (2008), Baltagi et al. (2012), and Judge and Mittelhammer (2012). When the instruments' exogeneity is in doubt, Conley et al. (2012) and Nevo and Rosen (2012) explore procedures to estimate bounds on the problematic regression coefficients.

This paper proposes the Non-Instrumental Simultaneous Equation (NISE) estimator, which requires no instruments and is consistent under a set of assumptions often invoked for linear statistical models. In an insightful paper on the various simultaneous-equation estimators, Chow (1964, 533-537) explains the relationship between canonical correlation (CC) and the estimator that I call NISE. He also compares NISE to LIML, noting that the latter has better statistical efficiency when valid instruments are available (Chow 1964, 542-543).

A rather long history of econometric research links CC to systems of linear equations. Early examples include Hooper (1959), who explores CC in relation to the reduced form, and Hannan (1967), who studies CC in the context of subsets of simultaneous equations that are just-identified. Hannan's paper inspired research by Chow and Ray-Chaudhuri (1967) and Ghosh (1972). Theil (1971, 317-321) uses CC to motivate Zellner's methodology for seemingly-unrelated regressions. More recently, CC and the closely-related reduced-rank regression (Reinsel and Velu 1998) have become the estimators for error-correction models of a vector autoregression that may contain cointegrated variables (e. g., Johansen 1996).

Section 2 introduces NISE as a maximum-likelihood estimator (mle) and describes a test for specification error; additional analytical content is provided in an appendix.  Section 3 reports a simulation of OLS, TSLS and NISE in the context of a market model, and section 4 applies the three estimation methods to a data set for Texas nursing facilities. Section 5 contains some conclusions and disclaimers.

**NISE**

The simultaneous linear equation of interest is

$$\mathbf{Y}\boldsymbol{\gamma} = \mathbf{X}\boldsymbol{\beta} + \mathbf{u}, \tag{1}$$

where the matrix **Y** contains n observations on G endogenous variables, **X** contains n observations on H exogenous variables, **γ** and **β** are vectors of unknown parameters, and **u** is a vector of n unobservable gaussian disturbances identically and independently distributed with $E(\mathbf{u}) = E(\mathbf{Xu}) = \mathbf{0}$. In addition to **X**, the model contains L exogenous variables **W** that appear in other linear equations; and L ≥ G so (1) is identifiable by exclusion restrictions.

A researcher wants to estimate (1) but does not intend to estimate the model's other equations. In fact, she may have no usable data on **W**, which is a principal motive for choosing NISE instead of LIML or TSLS. *It is irrelevant for identification that the actual measurements on the excluded exogenous variables may be absent from the data set or otherwise unusable.* If the model is correct, then the real effects of the excluded exogenous variables on the equation of interest are reflected in **Y**.

For this linear system, only one of whose equations is to be estimated, Davidson and MacKinnon (1993, 644) show that there is no Jacobian term in the loglikelihood; accordingly, the NISE mle is derived by minimizing the Lagrangian

$$F = (\mathbf{Y}\boldsymbol{\gamma} - \mathbf{X}\boldsymbol{\beta})^T(\mathbf{Y}\boldsymbol{\gamma} - \mathbf{X}\boldsymbol{\beta}) - \lambda[\boldsymbol{\gamma}^T(\mathbf{Y}^T\mathbf{Y})\boldsymbol{\gamma} - 1] , \tag{2}$$

where the constraint guarantees that the total squared error $\mathbf{u}^T\mathbf{u}$ is minimized with respect to the left-hand side of equation (1). (Chow 1964, 533-537 and 542-543; conforming to more recent usage, my notation reverses the roles of **β** and **γ** in his paper.) Chow shows that the estimate of **γ** is **c**, the characteristic vector corresponding to the smallest characteristic value $\lambda_{min}$ that solves the determinantal equation

$$|Y^TMY - \lambda Y^TY| = |Y^TMY(Y^TY)^{-1} - \lambda I| = 0. \qquad (3)$$

In (3), the nxn idempotent matrix $M = I-X(X^TX)^{-1}X^T$ "partials out" **X** from each column of **Y**.

At this point, a researcher may choose to normalize an endogenous variable with a coefficient of 1 on the equation's left-hand side, moving the other endogenous variable(s) to the right-hand side. For G ≥ 2 endogenous variables with $y_1$ on the left, the coefficient of $y_g$ is $-c_g/c_1$ (g = 2,…,G). Then **b**, an estimate of **β**, is obtained from the OLS regression of $y_1 + (c_2/c_1)y_2 + … + (c_G/c_1) y_G$ on **X**. Standard errors for the NISE coefficients can be estimated using a "pairs" or "cases" bootstrap.

The coefficient of determination from the OLS regression is $r_1^2 = 1 - \lambda_{min}$, the largest squared canonical correlation between **Y** and **X**. Accordingly, Bartlett's large-sample approximation to Wilks' lambda is a likelihood-ratio test for NISE:

$$Z = -(n-1-(G+H+1)/2)\sum \ln(1-r_i^2) \qquad , \qquad (4)$$

where the summation is over the sample squared canonical correlations $r_i^2$ from i = 2 to i = min(G,H). On the null hypothesis, only $r_1^2$ is statistically significant; that is, the data support exactly one linear relationship between **X** and **Y**, so equation (1) is correctly specified. Then Z has asymptotically a chi-square distribution with (G-1)(H-1) degrees of freedom (e. g., Mardia et al. 1979, 288-289).

A value of Z large enough to reject the null hypothesis at conventional levels of significance may be evidence of specification error and identification failure. If the null hypothesis in (4) is not rejected, NISE is consistent but inefficient since it is a limited-information method that ignores the model's other equations. NISE is also inefficient relative to LIML and TSLS methods when valid instruments exist. These inefficiency characterizations are of course based on asymptotic theory. Against them one can set the well-documented practical deficiencies of IV methods and the convenience of an additional reality check on IV estimates when they can be computed.

**A market simulation**

The simulation produces cross-section data on markets for an agricultural commodity. The demand function includes two endogenous variables, the price of the commodity (p) and the quantity demanded (qd); three exogenous variables –-household income (inc), the price of a substitute commodity (ps), and the price of a complementary commodity (pc); and a random disturbance $u_d$. The supply function includes the price and the quantity supplied (qs); the exogenous variables rainfall (r), the price of fertilizer (pf), and the ambient temperature (t); and a random disturbance $u_s$. The market is supposed to be in equilibrium at each observation in the data set: qd = qs = q. With the variables log-transformed, the simulated demand function is

$$q_i = -1.00 p_i + 1.50 inc_i + 0.50 ps_i - 0.50 pc_i + 3.00 + u_{di} \;; \quad (5)$$

and the simulated supply function is

$$q_i = 0.75 p_i + 2.50 r_i - 1.50 pf_i - 1.00 t_i + 0.50 + u_{si} \;, \quad (6)$$

for i = 1, …, n.

In both (5) and (6), the exogenous variables are uncorrelated with their respective disturbances; but the model implies that p and $u_d$ are indeed correlated. This is the cause of simultaneity bias that prevents consistent estimation by OLS. A researcher wants to estimate the demand function but not the supply function. However, the simulated values of the supply shifters r, pf and t will be required as instruments in the sequel since it is of interest to compare the performances of TSLS and NISE.

All the model's exogenous variables are generated from independent standard normal distributions; $u_d$ and $u_s$ are also independent normal, each with zero mean and a standard deviation of 2. At every replication, the simulation computes $p_i$ from the reduced form and inserts that value into equation (5). In 5,000 replications, the median correlation between p and $u_d$ is 0.489. Table 1 summarizes this base case for two sample sizes: n = 50 and n = 500. The table reports the median value of each regression coefficient and, underneath each coefficient, a robust version of its standard error, i. e., the Qn scale in Rousseeuw and Croux (1993).

    The simulations show that the OLS estimate of the price elasticity of demand is notably biased, and the other OLS coefficients are moderately biased toward zero. The TSLS and NISE coefficients are very close to the values in (5). For both sample sizes OLS is most efficient, followed by TSLS and then NISE; but when n = 500 the differences in efficiency are minor. The instrumental variables are strong as indicated by the significance level of the F test for the first-stage regression (it is a joint test on the exclusion of r, pf, and t). In addition, the null hypothesis for the J test is not rejected: r, pf and t have no significant correlation with the TSLS residuals. Finally, the null hypothesis for the Z test in NISE is not rejected: among the variables p, q, inc, ps and pc, the demand function is the only linear relationship supported by the data.

    Table 2 examines two variations on the base case. To explore the effect of *weak instruments* on TSLS, the coefficients of r, pf, and t in equation (6) are reduced to 0.5, -0.3, and -0.2 respectively. For this simulation, the median correlation between p and $u_d$ is 0.599; and at n = 50 the null hypothesis of the F test for weak instruments cannot be rejected. The TSLS price elasticity of demand is biased and statistically insignificant, but the NISE estimate is accurate and more than twice its robust standard error. At n = 500 the TSLS bias is much smaller and the standard error has shrunk.

    When the supply shifter *rain (r) is mistakenly included in the demand function*, the NISE price elasticity of demand is ruined at both sample sizes and its standard error explodes. This specification error is flagged by the Z test, for which the alternative hypothesis of non-identification cannot be rejected. TSLS performs well.

    In this section I have explicitly specified the demand and supply functions and have solved for the reduced form of the price variable, which I then substitute into the demand function to generate the quantity variable.

However, in their simulations of IV estimators some researchers (e.g., Staiger and Stock 1997, Judge and Mittelhammer 2012) do not actually specify the structural equations for the included endogenous variables and therefore do not derive an authentic reduced form. Instead the included endogenous variables are generated from a pseudo reduced form and then substituted into the equation of interest to produce the dependent variable on the left-hand side. Since these steps do not give rise to simultaneity bias, it is introduced *ad hoc* in the form of correlations between the structural disturbance and the disturbances in the pseudo reduced form.

This setup is effective for testing estimators like LIML and TSLS because they are general-purpose techniques, supposedly consistent in a wide variety of situations where a regressor is correlated with a disturbance. NISE, on the other hand, is not a general-purpose technique. It is focused on simultaneity bias, and for valid simulation it requires "data" from an authentic reduced form that is a function of the structural parameters and the structural disturbances.

**A cost function for nursing homes**

As an application of NISE, I draw on a data base of the Texas Health and Human Services Commission (2002) to estimate a short-run cost function for a typical Texas nursing home. For each of 920 facilities in 2002, the data set provides the total cost, the number of resident-days (the "quantity" variable), and the hourly wage rates paid to registered nurses (RN), licensed vocation nurses (LVN), nurse's aides (AIDE), and laundry and housekeeping personnel (L+H). The data were screened for outliers using robust methods (Maronna et al. 2006, 204-208), and 61 observations were removed.

Specifying a log-linear model, I assume that total cost and resident days are jointly endogenous since in theory management would try to minimize the cost per resident-day either as a primary objective or as a necessary condition for profit maximization. Because the nursing homes recruit staff in a large labor market of health-care personnel, I take the wage rates to be exogenous. The disturbances in the cost function are assumed to be independently and identically distributed normal variables, but they are correlated with resident days: given the wage rates, a nursing facility whose total cost is randomly higher or lower than expected will adjust its resident days.

Although NISE is likely to be most useful when there are no observations on the excluded exogenous variables, the data set appears to provide two candidates for instruments: the Medicaid per-diem reimbursement and the daily rate charged to private-pay patients. Microeconomic theory does not suggest that these variables belong in a cost function since the two payment rates pertain to the firm's revenue. The Medicaid per diem is arguably exogenous because it is set by regulatory authorities. However, the facility manager can affect the per diem by adjusting the composition of the resident days since the facility is reimbursed at higher rates for residents needing more assistance and supervision. The private-pay rate could also be exogenous when there is vigorous competition among nursing homes to enroll self-funded residents. On the other hand, the long-term care industry may resemble monopolistic competition if each facility can exercise some pricing power via the amenities it offers, the convenience of its location, its reputation for quality of care, and so on. In short, it is not certain that the candidate instruments are exogeneous.

Historically, occupancy rates in Texas nursing homes were lower than in many other states, partly because regulators readily authorized the construction of new or expanded facilities. The average occupancy rate in 2002 was just 70 percent, an indication that many firms were not operating very efficiently. The cost function could consequently be modeled using the stochastic-frontier methodology of Aigner et al. (1977) with its nonlinear likelihood, but I do not pursue this approach. As Greene (2003, 502-503) remarks, the asymmetry of the stochastic frontier's error distribution "does not negate our basic results for least squares in this classical regression model. This model satisfies the assumptions of the Gauss-Markov theorem…" apart from the above-mentioned endogeneity problem.

The cost-function estimates are shown in Table 3. As expected, all the slope coefficients are positive. The OLS and NISE coefficients are statistically significant (p-values < 0.10), but only the resident-days coefficient is significant for TSLS. The TSLS regression implies that total costs vary more than proportionally with resident days. For OLS, the relationship is essentially proportional; but total costs vary less than proportionally in the NISE regression. A bootstrap-based test indicates that each resident-days coefficient is significantly different from the other two.

In view of the previous discussion of the data set's instrumental variables, it is worth noting that the J test for the over-identifying restrictions is strongly rejected: the instruments are correlated with the TSLS residuals and, by inference, with the cost function's disturbances. So although the available instruments are strong --the F test in Table 3 clearly rejects the null hypothesis of weak instruments-- they don't verify the identification of the cost function. On the other hand, the null hypothesis of the Z test is definitely not rejected so the NISE estimates appear to be acceptable.

**Caveats and conclusions**

With the exception of cointegrated time series, for which OLS is "superconsistent," the econometrics literature appears to be united in the view that IV methods are the only way to estimate consistently a simultaneous linear equation when identification is based on exclusion restrictions. Nevertheless, this paper has provided tentative evidence that NISE can deal effectively with simultaneity bias when the linear equation is in fact identified and the sample size is not too small. Moreover, the Z test raises a red flag in situations where specification error causes a failure of identification. However, NISE cannot replace IV methods in applications other than simultaneous linear equations; for example, NISE can't deal with bias arising from "errors in variables." In addition, NISE fails when the equation of interest includes no exogenous variables: if **X** is empty, $\lambda_{min} = 1$ and also $\mathbf{Y}^T\mathbf{MY} = \mathbf{Y}^T\mathbf{Y}$ identically in equation (3). For the Z test, **X** must contain at least two variables.

For practical computations, the **X** and **Y** arrays can be input to CC software, which then produces **c** as the coefficients of **Y** corresponding to the *largest* squared canonical correlation. In addition, a pairs or cases bootstrap can estimate the sampling errors of the NISE coefficients, allowing for the possibility that the errors **u** in equation (1) are heteroskedastic. One can of course compute the usual bootstrap standard errors, but it may be preferable to use a robust dispersion estimator (Maronna et al. 2006, chapter 2, Rousseeuw and Croux 1993). In addition to limiting the distortions due to outlying data points, a robust bootstrap estimate of dispersion may be desirable since the NISE coefficients $-c_g/c_1$ do not necessarily have a finite second moment, as Anderson (2010) explains in the context of LIML.

**Appendix: a derivation of the NISE mle**

In equation (1), the covariance matrix of **Y**, denoted **Θ**, has rank G and is consistently estimated by its usual sample counterpart, denoted **V**. Likewise the covariance matrix of **X**, denoted **ψ**, has rank H. When the NISE log-likelihood is concentrated by "partialing out" **X** from each column of **Y,** the residuals **MY** have a covariance matrix $\sum$, of rank G, which is consistently estimated by its sample counterpart, denoted **S**. Then the variance of the jointly endogenous variables can be decomposed into the residual variance and the "regression" variance, say

$$\mathbf{\gamma}^T \mathbf{\Theta}\mathbf{\gamma} = \mathbf{\gamma}^T \textstyle\sum \mathbf{\gamma} + \mathbf{\delta}^T \mathbf{\psi}\mathbf{\delta} \ . \tag{A1}$$

Now NISE bears a formal resemblance to LIML, whose derivation is provided by Davidson and MacKinnon (1993, 644-649) and by Theil (1971, 502-504 and 679-686). LIML and NISE maximize similar constrained Gaussian log-likelihoods and have similar normalizations. For NISE, the normalization constraint is $\mathbf{\gamma}^T \mathbf{\Theta}\mathbf{\gamma} = 1$ (Chow 1964, 533-537 and 542-543; cp. Theil 1971, 684 and Anderson 2010, 359-361). However, NISE does not use the other LIML constraints, which involve the excluded exogenous variables **W**, i. e., the instruments. These latter constraints increase the asymptotic efficiency of LIML relative to NISE, but the premise of NISE is that valid observations on the instruments are not available.

In short, a derivation of NISE as mle follows the derivation of LIML if one omits every term involving variables excluded *a priori* from the equation of interest. Given a random sample of n observations **X** and **Y**, maximization of the NISE log likelihood is equivalent to the minimization of

$$F^* = 0.5n\log|\textstyle\sum| + 0.5\mathrm{tr}(\textstyle\sum^{-1}\mathbf{S}) - 0.5\rho(\mathbf{\gamma}^T\mathbf{\Theta}\mathbf{\gamma} - 1) \ , \tag{A2}$$

where ρ is a Lagrange multiplier (cp. Theil 1971, 679). **S** is a consistent estimator of ∑, and this substitution means that the second term in (A2) is a constant, not dependent on ∑. When (A1) is applied to the third term of (A2) and ∂F*/∂∑ is set equal to zero, it follows that

$$n\sum{}^{-1} - \rho \boldsymbol{\gamma}\boldsymbol{\gamma}^T = \mathbf{0} \qquad (A3)$$

(e.g., Theil 1971, 31-32 or Greene 2003, 839-840).

Premultiplying (A3) by ∑ and postmultiplying by **Θγ**,

$$\boldsymbol{\Theta}\boldsymbol{\gamma} - (\rho/n)\sum \boldsymbol{\gamma}\boldsymbol{\gamma}^T \boldsymbol{\Theta}\boldsymbol{\gamma} = \mathbf{0} \qquad (A4)$$

or $\quad [\boldsymbol{\Theta} - (\rho/n)\sum]\boldsymbol{\gamma} = \mathbf{0} \qquad (A5)$

or $\quad [\sum - \lambda\boldsymbol{\Theta}]\boldsymbol{\gamma} = \mathbf{0} , \qquad (A6)$

where λ = n/ρ. The corresponding determinantal equation is

$$|\sum - \lambda\boldsymbol{\Theta}| = |\sum\boldsymbol{\Theta}^{-1} - \lambda\mathbf{I}| = 0 ; \qquad (A7)$$

and minimization requires the smallest characteristic value,

$$\lambda_{min} = \boldsymbol{\gamma}^T \sum \boldsymbol{\gamma} / \boldsymbol{\gamma}^T \boldsymbol{\Theta}\boldsymbol{\gamma} , \qquad (A8)$$

which is a real positive fraction.

To operationalize the NISE estimation procedure, ∑ and **Θ** are replaced by their consistent estimators, **S** and **V** respectively. Then **γ** is estimated by **c**, the characteristic vector corresponding to $\lambda_{min}$. While standard errors for the NISE coefficients can be estimated using a bootstrap, it is also possible to estimate a covariance matrix for the coefficients based on $\tilde{\mathbf{y}}_g = \lambda_{min}\mathbf{y}_g - \mathbf{e}_g$, where $\mathbf{e}_g$ is the residual vector when **X** is partialed out of $\mathbf{y}_g$ (g = 2,…,G). If **A** is the sample covariance matrix of $(\tilde{\mathbf{y}}_2,...,\tilde{\mathbf{y}}_G,\mathbf{X})$, then an estimator of the asymptotic covariance matrix of $(-c_2/c_1,…,-c_G/c_1,\mathbf{b})$ is $n^{-1}\hat{\sigma}^2\mathbf{A}^{-1}$, where $\hat{\sigma}^2$ is the residual variance from the OLS regression that produces **b** as an estimate of **β** in section 2.

## Table 1. Base-case estimates of the demand function
## ( standard errors are under their coefficients)

|  | OLS | TSLS | NISE |
|---|---|---|---|
|  |  | n=50 |  |
| Median in simulation |  |  |  |
| Price | -0.603 | -0.973 | -1.004 |
|  | 0.109 | 0.165 | 0.234 |
| Income | 1.168 | 1.479 | 1.477 |
|  | 0.280 | 0.330 | 0.363 |
| Psub | 0.387 | 0.492 | 0.495 |
|  | 0.266 | 0.298 | 0.306 |
| Pcom | -0.381 | -0.487 | -0.486 |
|  | 0.268 | 0.301 | 0.308 |
| F signif |  | 0.000 |  |
| J signif |  | 0.494 |  |
| Z signif |  |  | 0.505 |
|  |  | n = 500 |  |
| Price | -0.600 | -0.999 | -0.998 |
|  | 0.033 | 0.050 | 0.068 |
| Income | 1.157 | 1.496 | 1.496 |
|  | 0.084 | 0.100 | 0.106 |
| Psub | 0.385 | 0.502 | 0.501 |
|  | 0.079 | 0.091 | 0.093 |
| Pcom | -0.386 | -0.499 | -0.498 |
|  | 0.079 | 0.091 | 0.091 |
| F signif |  | 0.000 |  |
| J signif |  | 0.489 |  |
| Z signif |  |  | 0.499 |

**Table 2. Additional estimates of the price elasticity of demand**
**( standard errors are under their coefficients)**

|  | OLS | TSLS | NISE |
|---|---|---|---|
| Median in simulation | | Weak instruments, n = 50 | |
| elasticity | -0.167 | -0.520 | -0.999 |
|  | 0.130 | 0.531 | 0.416 |
| F signif |  | 0.323 |  |
| J signif |  | 0.499 |  |
| Z signif |  |  | 0.502 |
|  | | Weak instruments, n = 500 | |
| elasticity | -0.164 | -0.927 | -1.002 |
|  | 0.040 | 0.232 | 0.118 |
| F signif |  | 0.003 |  |
| J signif |  | 0.499 |  |
| Z signif |  |  | 0.496 |
| Median in simulation | | RAIN in demand function, n = 50 | |
| elasticity | -0.378 | -0.946 | -0.031 |
|  | 0.126 | 0.278 | 2.217 |
| F signif |  | 0.000 |  |
| J signif |  | 0.499 |  |
| Z signif |  |  | 0.000 |
|  | | RAIN in demand function, n = 500 | |
| elasticity | -0.378 | -0.999 | 0.534 |
|  | 0.038 | 0.088 | 1.074 |
| F signif |  | 0.000 |  |
| J signif |  | 0.500 |  |
| Z signif |  |  | 0.000 |

**Table 3. Estimates of the cost-function model**
**the dependent variable is log total cost**
**(standard errors are under their coefficients)**
**n = 859**

|  | OLS | TSLS | NISE |
|---|---|---|---|
| ln resident days | 0.977 | 1.216 | 0.779 |
|  | 0.009 | 0.079 | 0.025 |
| ln RN wage | 0.083 | 0.074 | 0.090 |
|  | 0.042 | 0.057 | 0.055 |
| ln LVN wage | 0.271 | 0.033 | 0.468 |
|  | 0.050 | 0.102 | 0.070 |
| ln aide wage | 0.261 | 0.116 | 0.381 |
|  | 0.047 | 0.079 | 0.058 |
| ln L+H wage | 0.098 | 0.055 | 0.133 |
|  | 0.053 | 0.072 | 0.067 |
| F signif 1st-stage reg |  | 0.000 |  |
| J signif |  | 0.000 |  |
| Z signif |  |  | 0.592 |

|  | OLS -NISE | TSLS-OLS | TSLS-NISE |
|---|---|---|---|
| ln resident days | 0.198 | 0.239 | 0.437 |
|  | 0.024 | 0.128 | 0.131 |